\documentclass[a4paper,12pt]{article}
\usepackage[top=2.5cm,bottom=2.5cm,left=2.5cm,right=2.5cm]{geometry}
\usepackage{cite, amsmath, amssymb}
\usepackage[margin=1cm,%
            font=small,%
            format=hang,%
            labelsep=period,%
            labelfont=bf]{caption}
\usepackage{graphicx}
\usepackage{cite}
\allowdisplaybreaks[1]

\begin{document}
\title{\vspace{3.5 cm} On the Maximum ABC Spectral Radius\\of Connected Graphs and Trees\footnote{Supported by the National Science Foundation of China under Grant Nos. 11771362 and 11601006.} }

\author{Wenshui Lin$^{1,2}$, Yiming Zheng$^2$, Peifang Fu$^2$, Zhangyong Yan$^2$, Jia-Bao Liu $^{3,}$\footnote{Corresponding author. E-mail address: liujiabaoad@163.com (J. Liu).}\\
        {\small $^1$ Fujian Key Laboratory of Sensing and Computing for Smart City, Xiamen 361005, China}\\
        {\small $^2$ School of Informatics, Xiamen University, Xiamen 361005, China}\\
        {\small $^3$ School of Mathematics and Physics, Anhui Jianzhu University, Hefei 230601, China}
}

\date{\small (Received April 16, 2020)}
\maketitle \thispagestyle{empty}

\noindent \textbf{Abstract}\\
Let $G=(V,E)$ be a connected graph, where $V=\{v_1, v_2, \cdots, v_n\}$ and $m=|E|$. $d_i$ will denote the degree of vertex $v_i$ of $G$, and $\Delta=\max_{1\leq i \leq n} d_i$. The ABC matrix of $G$ is defined as $M(G)=(m_{ij})_{n \times n}$, where $m_{ij}=\sqrt{(d_i + d_j -2)/(d_i d_j)}$ if $v_i v_j \in E$, and 0 otherwise. The largest eigenvalue of $M(G)$ is called the ABC spectral radius of $G$, denoted by $\rho_{ABC}(G)$. Recently, this graph invariant has attracted some attentions. We prove that $\rho_{ABC}(G) \leq \sqrt{\Delta+(2m-n+1)/\Delta -2}$. As an application, the unique tree with $n \geq 4$ vertices having second largest ABC spectral radius is determined.

\baselineskip=0.30in

\noindent \textbf{Keywords:} ABC matrix, Eigenvalues, ABC spectral radius, Upper bounds, Trees.

\baselineskip=0.30in

\section{Introduction}
Let $G=(V,E)$ be a simple connected graph. Suppose $V=\{v_1, v_2, \cdots, v_n\}$ and $m=|E|$. If $c=m-n+1$ ($\geq 0$), then $G$ is called a $c$-cyclic graph. In particular, $G$ is called a tree and a unicyclic graph if $c=1, 2$, respectively. As usual, $S_n$, $P_n$, $C_n$, and $K_n$ will denote the star, path, cycle, and complete graph with $n$ vertices, respectively.

Let $d_i$ denote the degree of vertex $v_i$, and $\Delta=\max_{1\leq i \leq n} d_i$.
The atom-bond connectivity index (ABC index in short) of $G$ is defined \cite{[1]} as
$ABC(G) = \Sigma_{v_i v_j \in E} f(d_i, d_j)$, where $f(x,y)=\sqrt{ (x+y-2)/(xy) }$.
Since this index can predict well the heat of formation of alkanes (see \cite{[2],[3]}), it became a hot topic in the past few years (see [4-30]).

In 2017, Estrada \cite{[31]} defined the ABC matrix of $G$ as $M = M(G)=(m_{ij})_{n \times n}$, where $m_{ij}= f(d_i, d_j)$ if $v_i v_j \in E$, and $0$ otherwise. The chemical background of this matrix was explicated in \cite{[31]}.
The eigenvalues of $M$ are called the ABC eigenvalues of $G$. Because $M$ is non-negative, symmetric, and irreducible, any ABC eigenvalue of $G$ is real.
In particular, the largest ABC eigenvalue of $G$ is called its ABC spectral radius,
and denoted by $\rho_{ABC} (G)$. Obviously, $\rho_{ABC} (G)$ is positive and simple.
Moreover, there exists a unique vector $x >0$ such that $\rho_{ABC} (G) = \max_{\|y\|=1} y^T My = x^T Mx$,
which is known as the Perron vector of $M$.

Estrada \cite{[31]} proved that $\frac{2}{n} ABC(G) \leq \rho_{ABC} (G) \leq \max_{1 \leq i \leq n} M_i$,
with both equalities iff $M_1 = M_2 = \cdots = M_n$, where $M_i = \sum_{1 \leq j \leq n} m_{ij}$.
Recently, Chen \cite{[32]} presented another lower bound of $\rho_{ABC} (G)$ in terms of $R_{-1} (G)$,
which is the sum of $\frac{1}{d_i d_j}$ over all edges $v_i v_j \in E$.
Chen \cite{[32]} further proposed the problem of characterizing graphs
with extremal ABC spectral radius for a given graph class.
Soon, this problem for trees, connected graphs, and unicyclic graphs were solved by
Chen \cite{[33]}, Ghorbani et al. \cite{[34]}, and Li et al. \cite{[35]}, respectively.

\noindent \textbf{Lemma 1.1} \cite{[33]}\textbf{.} Let $T$ be a tree with $n \geq 3$ vertices. Then
$$ \sqrt{2} \cos\frac{\pi}{n+1} \leq \rho_{ABC} (G) \leq \sqrt{2n-4},$$
\noindent with left (right) equality iff $G \cong P_n$ (resp. $G \cong S_n$).

\noindent \textbf{Lemma 1.2} \cite{[34]} \textbf{.} Let $G$ be a connected graph with $n \geq 3$ vertices. Then
$$\sqrt{2} \cos{ \frac{\pi}{n+1} \leq \rho_{ABC} (G) \leq \sqrt{ 2n-4 }},$$
with the left (right) equality iff $G \cong P_n$ (resp. $G \cong K_n$ ).

\noindent \textbf{Lemma 1.3} \cite{[35]} \textbf{.} Let $G$ be a unicyclic graph with $n \geq 4$ vertices. Then
$$ \sqrt{2} \ = \rho_{ABC} (C_n) \leq \rho_{ABC} (G) \leq \rho_{ABC} (S_n +e),$$
with the left (right) equality iff $G \cong C_n$ (resp. $G \cong S_n +e$).

For convenience, let $\mathcal{C}(n)$ be the set of connected graphs with $n$ vertices,
and $\mathcal{G}(m,n)$ the set of connected graphs with $n$ vertices and $m$ edges.
In the present paper, we consider upper bounds of $\rho_{ABC}$ for connected graphs.
In Section 2, it is shown that, if $G \in \mathcal{G}(m,n)$ and $\Delta(G) = \Delta$,
then $\rho_{ABC}(G) \leq \sqrt{\Delta+(2m-n+1) / \Delta -2}$.
As an application, in Section 3, we characterize the unique tree with $n \geq 4$ vertices
having the second largest ABC spectral radius. Finally, some problems are proposed in Section 4.

\section{Some upper bounds of the ABC spectral radius}

In this section, we present two upper bounds of $\rho_{ABC}$ of connected graphs.

\noindent \textbf{Theorem 2.1.} If $G \in \mathcal{G}(m,n)$ and $\Delta(G) = \Delta$, then
$$\rho_{ABC} (G) \leq \sqrt{\Delta + (2m-n+1)/ \Delta -2}.$$
\noindent Moreover, the bound is attainable.

\noindent \textbf{Proof.} Let $M = M(G)$, $D = 2m -n +1$, and $x = ( \sqrt{d_1}, \sqrt{d_2}, \cdots, \sqrt{d_n} )^T $. From the Perron-Frobenius theory, it suffices to confirm the claim: $(Mx)_i \leq \sqrt{ d_i} \sqrt{ \Delta + D / \Delta -2 }$ or $[ (Mx)_i / \sqrt{d_i}] ^2 \leq \Delta + D/ \Delta -2 $ holds for $ 1 \leq i \leq n$.

If $d_i < D/ \Delta$, then $ (Mx)_i = \sum \limits_{v_i v_j \in E} f(d_i, d_j) \sqrt{d_j}
\leq d_i \sqrt{ \frac{d_i +\Delta -2}{d_i} }
< \sqrt{d_i} \sqrt{ D / \Delta + \Delta-2}$. Hence assume $d_i \geq D/ \Delta$.
By using the Cauchy-Schwarz Inequality we have
\begin{align*}
  (Mx)_i &= \sum \limits_{v_i v_j \in E} \sqrt{ (d_i + d_j -2) / d_i } \\
         &\leq \sqrt{ d_i^2 -2d_i + \sum \limits_{v_i v_j \in E} d_j } \\
         &\leq \sqrt{ d_i^2 -2d_i + [2m-d_i-(n-d_i-1)] } \\
         &= \sqrt{ d_i^2 -2d_i + D }.
\end{align*}

\noindent Thus we have $[ (Mx)_i / \sqrt{d_i}] ^2 \leq (d_i^2 -2d_i +D)/d_i = d_i + D/d_i -2 $.
Since $\eta (x) = x+D/x -2$ is a Nike function and $D/ \Delta \leq d_i \leq \Delta$, it follows that
$$[ (Mx)_i / \sqrt{d_i}] ^2 \leq \max \{\eta (D/ \Delta), \eta (\Delta) \} = \Delta + D/ \Delta -2.$$

Finally, to see the bound is attainable, one can take $S_n$ and $K_n$ as examples. The proof is thus completed. $\blacksquare$

Let $\theta(m, \Delta) = \sqrt{ \Delta +(2m-n+1)/ \Delta -2 }$. For fixed $m$, the monotonicity of $\theta$ with respect to $\Delta$ is clear. Hence Theorem 2.1 can easily produce a upper bound of $\rho_{ABC}$ for subsets of $\mathcal{G}(m, n)$. For example, an upper bound for $c$-cyclic graphs is obtained as follows.

\noindent \textbf{Corollary 2.2.} Let $G$ be a $c$-cyclic graph with $n \geq 3$ vertices, and $c \leq (n-1)/2$. Then
$$\rho_{ABC} (G) \leq \sqrt{n-2 + 2c/(n-1)}.$$

\noindent \textbf{Proof.} Since $m=n-1+c$ and $c \leq (n-1)/2$, by direct calculations we have
$$\theta (m, 2) = \sqrt{ (n-1)/2 +c } \leq \theta (m, n-1) = \sqrt{ n-2 +2c/(n-1) },$$

\noindent and the conclusion follows from Theorem 2.1. $\blacksquare$

It is easily seen that, Theorem 2.1 can reproduce the upper bound part of Lemma 1.1. However, if we consider upper bounds of $\rho_{ABC}$ for a subset of $\mathcal{C}_n$, whose elements have various sizes (numbers of edges),
Theorem 2.1 may be not so convenient to applied directly. Hence we deduce the following result.

\noindent \textbf{Corollary 2.3.} If $G \in \mathcal{G}(m,n)$ and $\Delta(G) = \Delta$, then $\rho_{ABC} (G) \leq \sqrt{ \Delta +k-2 }$, where $k= \lceil (2m-n+1)/ \Delta \rceil$.

\noindent \textbf{Proof.} From $k= \lceil (2m-n+1)/ \Delta \rceil \geq (2m-n+1)/ \Delta$,
the conclusion holds immediately from Theorem 2.1. $\blacksquare$

Though Corollary 2.3 is weaker than Theorem 2.1, the upper bound $\theta'(m, \Delta) = \sqrt{ \Delta +k-2 }$  has better property than $\theta(m, \Delta)$. In fact, for fixed $m$, $\theta'(m, \Delta)$ almost strictly decreases with $\Delta$. To see this, let $\Delta_1 > \Delta_2$ and $D = 2m-n+1$. We illustrate the fact with the following two cases.

\textbf{Case 1.} $D/k \leq \Delta_2 < \Delta_1 < D/(k-1)$. Then $\theta'(m, \Delta_2) < \theta'(m, \Delta_1)$ obviously.

\textbf{Case 2.} $D/k \leq \Delta_2 < D/(k-1) \leq \Delta_1$. Then
\begin{align*}
  \theta'(m, \Delta_2)  &= \sqrt{\Delta_2 +k-2}\\
                        &\leq \sqrt{ \Delta_1 +k-3 }\\
                        & = \theta'(m, \Delta_1),
\end{align*}
\noindent with equality iff $\Delta_2 = D/(k-1) -1 = \Delta_1 -1$.

By the monotonicity of $\theta'(m, \Delta)$ with respect to $\Delta$, we are able to reproduce the upper bound part of Lemma 1.2.

\noindent \textbf{Corollary 2.4} \cite{[34]} \textbf{.} Let $G$ be a connected graph with $n \geq 3$ vertices. Then
$$\rho_{ABC}(G) \leq \sqrt{2n-4},$$
\noindent with equality iff $G \cong K_n$.

\noindent \textbf{Proof.} By the monotonicity of $\theta'$ we have $\theta'( m, \Delta ) \leq \theta'( m, n-1 ) \leq \theta'( n(n-1)/2, n-1 ) $, with all the equalities iff $m= n(n-1)/2$, that is, $G \cong K_n$. The conclusion thus follows from Corollary 2.3. $\blacksquare$

\section{The tree with second largest ABC spectral radius}

In order to further illustrate the application of Theorem 2.1 and Corollary 2.3, in this section we determine the tree with $n \geq 4$ vertices, whose ABC spectral radius is the second largest. For convenience, let $\mathcal{T}_n$ be the set of trees with $n$ vertices, and $\mathcal{T}_n ^{ (\Delta)}= \{ T \in \mathcal{T}_n| \Delta(T) = \Delta \}$.
Let $T_i$, $i=1, 2, 3, 4$, be the trees shown in Figure 1. $T_1$ is just the double star $S_{n-3, 1}$.
If $T \in \mathcal{T}_n ^{ (\Delta) }$, then $T$ contains $S_{\Delta +1}$ as its (induced) subgraph,
hence it is easily seen that $\mathcal{T}_n ^{ (n-1)} = \{S_n \}$, $\mathcal{T}_n ^{ (n-2)} = \{S_{n-3, 1}\}$,
and $\mathcal{T}_n ^{ (n-3)} = \{T_2, T_3, T_4 \}$.

\begin{figure}
  \centering
  \includegraphics{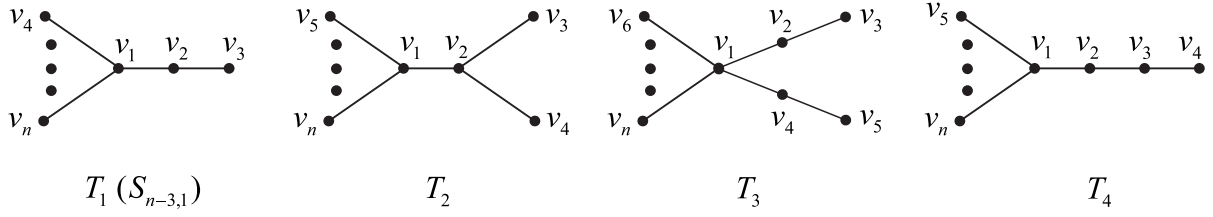}
  \caption{The trees $T_i$, i=1, 2, 3, 4.}
\end{figure}

Our aim in this section is to prove the following conclusion.

\noindent \textbf{Theorem 3.1.} If $n \geq 4$ and $T \in \mathcal{T}_n - \{ S_n, S_{n-3, 1} \}$, then
$$\rho_{ABC}(T) <  \rho_{ABC}(S_{n-3, 1}) < \rho_{ABC}(S_n).$$

We need some more preliminaries before presenting the proof of Theorem 3.1.

For two vertices $u$ and $v$ of a graph $G$, they are said to be equivalent, denoted by $u \sim v$,
if there is an automorphism of $G$ sending $u$ to $v$. By symmetry, the following result is immediate.

\noindent \textbf{Lemma 3.2.} Let $x$ be the Perron vector of the ABC matrix $M(G)$ of a connected graph $G$.
If $u \sim v$, then $x_i = x_j$.

\noindent \textbf{Lemma 3.3.} $\rho_{ABC}(S_{n-3, 1}) > \sqrt{n-3.5}$ if $n \geq 4$.

\noindent \textbf{Proof.} Let $\rho = \rho_{ABC}(S_{n-3, 1})$, and label the vertices of $S_{n-3, 1}$ as in Figure 1. Based on Lemma 3.2, let $x=(x_1, x_2, x_3, x_4, \cdots, x_4)^T$ be the Perron vector of $M = M(S_{n-3, 1})$. From $\rho x = Mx$ we have

\begin{equation*}
\left\{
  \begin{aligned}
    &\rho x_1 = (n-3) \sqrt{ \frac{n-3}{n-2} }x_4 + \sqrt{ \frac{1}{2}}x_2\\
    &\rho x_2 = \sqrt{ \frac{1}{2}}x_1 + \sqrt{ \frac{1}{2}}x_3 > \sqrt{ \frac{1}{2}}x_1\\
    &\rho x_4 = \sqrt{ \frac{n-3}{n-2} }x_1
  \end{aligned}
\right..
\end{equation*}

\noindent Hence $\rho^2 x_1 = (n-3) \sqrt{ \frac{n-3}{n-2} } \rho x_4 + \sqrt{ \frac{1}{2} } \rho x_2
> \frac{ (n-3)^2 }{n-2} x_1 + \frac{1}{2} x_1$, and we arrive at
$$\rho^2 > \frac{ (n-3)^2 }{n-2} + \frac{1}{2}
> n-4+0.5 = n-3.5,$$

\noindent which completes the proof. $\blacksquare$

\noindent \textbf{Lemma 3.4.} If $n \geq 6$ and $T \in \{T_2, T_3, T_4 \}$,
then $\rho_{ABC} (T) < \sqrt{n-3.5}$.

\noindent \textbf{Proof.} If $n=6$, the conclusion can be verified easily. Hence assume $n \geq 7$. Label the vertices of $T$ as in Figure 1. Let $M=M(T)$ and $x=( \sqrt{d_1}, \sqrt{d_2}, \cdots, \sqrt{d_n} )^T$.
Based on Lemma 2.1, we prove the result by confirming $(Mx)_i / \sqrt{d_i} \leq \sqrt{ n-3.5 }$ for $1 \leq i \leq n$.

We have $(Mx)_i / \sqrt{d_i} \leq \sqrt{ d_i-2+ \sum_{v_i v_j \in E} d_j/d_i}$ from the proof of Theorem 2.1. Hence $(Mx)_1 / \sqrt{d_1} \leq \sqrt{ n-5 +(n-1)/(n-3)} = \sqrt{ n-4 +2/(n-3) } \leq \sqrt{ n-3.5}$. For $i \geq 2$, because $d_i \leq 3$ and $n \geq 7$, it follows
$$(Mx)_i / \sqrt{d_i} \leq \max\{ \sqrt{n-4}, \sqrt{(n-1)/2}, \sqrt{(n+2)/3} \} < \sqrt{n-3.5},$$

\noindent and the proof is completed. $\blacksquare$

Now we present the proof of Theorem 3.1.

\noindent \textbf{Proof of Theorem 3.1.} It is easily seen that $\mathcal{T}_4 = \{ S_4, S_{1, 1} \cong P_4 \}$,
$\mathcal{T}_5 = \{ S_5, S_{2, 1}, P_5 \}$, and $\mathcal{T}_6 = \{ S_6, S_{3, 1}, T_2, T_3, T_4, P_6 \}$,
so the conclusion holds if $n \leq 6$ from Lemmas 1.1, 3.3, and 3.4. Hence assume $n \geq 7$.

If $\Delta \geq n-3$, the conclusion follows from Lemmas 3.3 and 3.4. Otherwise, if $\Delta \leq n-4$, then $\rho_{ABC}(T) \leq \theta' (n-1, n-4) = \sqrt{n-6+ \lceil (n-1)/(n-4) \rceil } \leq \sqrt{n-4}$ from Corollary 2.3.

The proof is thus completed. $\blacksquare$

\section{Further discussions}
In this paper, it is shown that $\rho_{ABC}(G) \leq \sqrt{\Delta + (2m-n+1)/ \Delta -2}$ if $G \in \mathcal{G}(m, n)$ and $\Delta(G) = \Delta$. The bound is attained by $S_n$ and $K_n$. Firstly, the following problem may be worth consideration.

\noindent \textbf{Problem 4.1.} Characterize the graphs $G \in \mathcal{G}(m,n)$ such that
$$\rho_{ABC}(G) = \sqrt{\Delta(G) + (2m-n+1)/ \Delta(G) -2}.$$

As well, the double star $S_{n-3, 1}$ is shown to be the unique tree having second largest ABC spectral radius in $\mathcal{T}_n$, $n \geq 4$. Recall that, Lin et al. \cite{[36]} ordered trees by their (adjacent) spectral radius $\lambda_1$, and showed that, if $T_1$ and $T_2$ are two trees with $n \geq 4$ vertice and $\Delta(T_1) > \Delta(T_2) \geq (2n)/3-1$, then $\lambda_1(T_1) > \lambda_1(T_2)$.
Naturally, the following question is interesting.

\noindent \textbf{Question 4.2.} Let $G_1$ and $G_2$ be two graphs in a subset of $\mathcal{G}(m, n)$. Is there some integer $l(m, n)$ (depending on $n$ and/or $m$), such that if $\Delta(G_1) > \Delta(G_2) \geq l(m,n )$, then $\rho_{ABC}(G_1) > \rho_{ABC}(G_2)$?

This question may be difficult to answer at the present, even for trees, and the following two problems are worth investigation in advance.

\noindent \textbf{Problem 4.3} Order graphs in some classes of connected graphs by their ABC spectral radii.

\noindent \textbf{Problem 4.4.} Establish non-trivial lower bounds of $\rho_{ABC}(G)$ for a graph $G \in \mathcal{G}(m, n)$ (in terms of $m$, $n$, and $\Delta(G))$.

\end{document}